\documentclass[a4paper,11pt]{article}
\usepackage{authblk}
\usepackage[latin1]{inputenc}
\usepackage[pdftex]{graphicx}
\usepackage{subfig}
\usepackage{float}
\usepackage{amsfonts,amssymb,amsmath,amsthm}
\usepackage{pdfsync}
\usepackage{color}
\usepackage{url}
\usepackage{nicefrac}
\usepackage[pdftex]{hyperref}
\usepackage{hyperref}

 \hypersetup{colorlinks=true,citecolor=blue, filecolor=red, linkcolor=blue, urlcolor=blue}

\newtheorem{theorem}{Theorem}[section] 
\newtheorem{proposition}[theorem]{Proposition}

\newtheorem{definition}[theorem]{Definition}
\newtheorem{remark}[theorem]{Remark}

\makeatletter
\@addtoreset{equation}{section}

\makeatother

\numberwithin{equation}{section}

\newcommand{\R}{\mathbb{R}}
\newcommand{\N}{\mathbb{N}}
\newcommand{\Z}{\mathbb{Z}}

\def\ds{\displaystyle}

\def\rg{\rightarrow}
\def\e{\varepsilon}

\newcommand{\be}{\begin{equation}}
\newcommand{\ee}{\end{equation}}
\newcommand{\ben}{\begin{equation*}}
\newcommand{\een}{\end{equation*}}
\newcommand{\ba}{\begin{eqnarray}}
\newcommand{\ea}{\end{eqnarray}}
\newcommand{\ban}{\begin{eqnarray*}}
\newcommand{\ean}{\end{eqnarray*}}

 \makeindex

\begin{document}

\title{Quasimodes and unstability for linear Schr\"odinger equation on manifolds} 
\author{Philippe Kerdelhué}
 \affil{ \footnotesize   Département de Mathématiques, CNRS UMR 8628, \\ F-91405 Orsay Cedex, France}

\date{}

\maketitle

  \begin{abstract} 
We consider the evolution operator $\exp(-it(-\Delta+V))$ associated with a Schr\"odinger operator on a Riemannian manifold $(M,g)$.
We are interested in the dependence of this operator on $V$ running in $L^p(M)$. Under a geometrical hypothetis,
we show the unstability for $p<\infty$ and give examples for which the hypothetis is satisfied. Then we show in the
general case the unstability for $p<\dim M/2$.
 \end{abstract}
 
\section[Introduction]{Introduction} \label{intro} 

These last years several papers were published about the stability of the nonlinear Schr\"odinger equation,
that is the uniform continuity of the map $u_0\mapsto u$ where $u$ is the solution of the equation
\be i\partial_t u+\Delta u=\varepsilon |u^2|u,~~~~u(0,x)=u_0 \label{nls} \ee
 where $\varepsilon=\pm 1$, $u(t,\cdot)$ is defined
on a Riemannian manifold $M$ and $\Delta$ is the Laplace-Beltrami operator on this manifold. More precisely :
\begin{definition} \label{defnls}
 Let $\sigma\in\R$, and denote by $B_{R,\sigma}$ the ball of radius $R$ in $H^\sigma(M)$. 
We say that the problem \ref{nls} is uniformly well-posed in $H^\sigma$ if for any $R>0$,
there exists $T>0$ suth that the map :
\ben B_{R,\sigma}\cap H^1(M)\ni u_0\mapsto u\in L^\infty([-T,T];\,H^\sigma (M)) \een
in uniformly continuous ($B_{R,\sigma}\cap H^1(M)$ is endowed with the $H^\sigma$ norm).\\
Otherwise, we say that the Cauchy problem (\ref{nls}) is unstable.
\end{definition}
Let us recall some results :\\
In 1993, J.Bourgain proved in \cite{Bo} that the Cauchy problem is uniformly well-posed
on the rational torus $\mathbb T^2$ when $\sigma>0$,\\
in 2002, N.Burq, P.G\'erard and N.Tzevtkov showed in \cite{BGT1} the unstability on $\mathbb S^2$ when $\ds 0\leq\sigma<1/4$,\\
in 2004, the same authors proved in \cite{BGT2} that the Cauchy problem is uniformly well-posed when
 $M$ is a compact manifold of
dimension $d\geq 2$ and $\ds \sigma>(d-1)/2$,\\
in 2005, they showed in \cite{BGT3} that the Cauchy problem is well-posed in $\mathbb S^2$ when $\sigma>1/4$,\\
in 2008, L.Thomann \cite{Tho} proved the unstability in the case when $M$ is a surface with a stable and not degenerate periodic
geodesic and $0<\sigma<1/4$.\\
For results in one dimension, we refer to the 
work of M.Christ, J.Colliander and T.Tao \cite{CCT}.\\

In the linear case, the propagators are unitary operators so the dependence on the initial data has no interest.
However, the stability property is the continuity of the evolution operator $\exp({-it(-\Delta+V)})$ with respect to $V$ running
in a proper functions space. Recently in \cite{BBZ} J.Bourgain, N. Burq and M.Zworski proved the following stability
result on the torus $\mathbb T^2=\R^2/(a\Z\times b\Z)$, $(a,b)\in\R^2$

\begin{theorem} \label{thBBZ}
\begin{enumerate}
 \item Let $K$ be a compact subset of $L^2(\mathbb T^2)$. Then for any $T>0$ the map
\ben K\ni V \mapsto e^{-it(-\Delta+V)} \in \L^\infty((0,T),\mathcal L (L^2(\mathbb T))) \een
is lipschitz continuous.\\
\item Let $p>2$ and $A$ be a bounded subset of $L^p(\mathbb T^2)$. Then for any $T>0$ the map
\ben A\ni V \mapsto e^{-it(-\Delta+V)} \in \L^\infty((0,T),\mathcal L (L^2(\mathbb T))) \een
is lipschitz continuous.
\end{enumerate}
\end{theorem}

The authors notice that ``it would be interesting to investigate such properties on other manifolds, as they seem to depend
strongly on the geometry''. The aim of this article is to answer partially to this remark.\\

This article is organized as following

\begin{itemize}
\item In section \ref{stat} we state our main theorem : assuming a geometrical condition on a sequence
of quasimodes for a Schr\"odinger operator $-\Delta+V$ on a manifold $M$, we show the unstability near $V$ of the maps
$ W \mapsto e^{-it(-\Delta+W)} \in \L^\infty((0,T),\mathcal L (L^2(\mathbb T)))$ with respect to the $L^p$ norm, $1\leq p<+\infty$.
\item In section \ref{exa} we show examples for which this theorem applies.
\item In section \ref{ppetit} we prove that the geometrical hypothetis is not necessary for $\ds p<\dim M/2$.
\end{itemize}

{\bf Acknoledgements :} The author would like to thank N.Burq for suggesting this subject and for useful help in the realization of
this article.

\section{Main theorem} \label{stat}

\begin{theorem} \label{mainth}

 Let $(M,g)$ be a Riemannian manifold of finite dimension, $\Delta_g=\Delta$ the Laplace-Beltrami operator, $d\mu$ the riemannian
volume form and $V$ a continuous nonnegative potentiel.
Assume
\begin{enumerate} 
\item  There exists a sequence $(\lambda_n,u_n)$ of quasi-eigenvalues and associated quasimodes for the Schr\"odinger operator
$-\Delta +V$ 
$$\lim_{n\rg +\infty} \|(-\Delta+V-\lambda_n)u_n\|_{L^2(M)}=0,~~\|u_n\|_{L^2(M)}=1.$$
\item The sequence of measures $|u_n(x)|^2 d\mu$ tends to a measure $\nu$ for the weak* topology, that is for any
continuous $f$ vanishing at infinity
$$\lim_{n\rg +\infty} \int_M f(x) |u_n(x)|^2 d\mu= \int_M f(x) \, d\nu$$
and $\nu$ is not absolutely continuous with respect to $\mu$.
\end{enumerate}
Then there exists a sequence of smooth bounded potentiels $W_k$ such that 
$$\|W_k\|_{L^{\infty}(M)}\leq 1$$
$$\forall p\in [1,+\infty[,~~\lim_{k\rg +\infty} \|W_k\|_{L^{p}(M)}=0$$
$$\forall t>0,~~\liminf_{k\rg +\infty} \|e^{-it(-\Delta+V+W_k)}- e^{-it(-\Delta+V)}\|_{\mathcal L(L^2(M))} >0$$

\end{theorem}

\vskip 5mm

{\bf Proof :} First notice that the operators $-\Delta+V$ and $-\Delta+V+W_n$ are self-adjoined with the same
domain $\{u\in\L^2(M);\,\Delta\,u\in L^2(M),\,V\,u\in L^2(M)\}$.\\
The second hypothetis ensures the existence of a compact subset $\Gamma$ of $M$ such that
$\mu(\Gamma)=0$ and $\nu(\Gamma)>0$. Let $\kappa\in]0,1]$ be a parameter to be fixed later and
$\left(\varphi_k\right)_{k\in\N^*}$ be a sequence of smooth fonctions on $M$ with values in $[0,1]$ such that
$$\left\{\begin{array}{l}
\varphi_k=1\text{ on }\Gamma\\
\mu(\text{supp}(\varphi_k))\leq\frac{1}{k}
\end{array}\right.$$
We observe 
$$\displaylines{
(-i\partial_t+(-\Delta+V(x)+\kappa\varphi_k(x)))\left(e^{-i(\lambda_n+\kappa)t} u_n(x)\right) \hfill \cr
\hfill =e^{-i(\lambda_n+\kappa)t}(-\Delta+V(x)-\lambda_n) u_n(x)+\kappa(\varphi_k(x)-1)\,e^{-i(\lambda_n+\kappa)t} u_n(x)\cr}$$
By hypothetis 1
$$\lim_{n\rg +\infty} \|e^{-i(\lambda_n+\kappa)t}(-\Delta+V-\lambda_n) u_n\|_{L^2(M)}=0$$
On the other hand 
$$\|(\varphi_k(x)-1)e^{-i(\lambda_n+\kappa)t} u_n(x)\|^2_{L^2(M)}=1-\int_M (2\varphi_k(x)-\varphi_k(x)^2)|u_n(x)|^2 d\mu$$
tends to $\ds 1-\int_M (2\varphi_k-\varphi_k^2)\, d\nu$ when $n$ tends to $+\infty$, and this quantity tends to $1-\nu(\Gamma)$
when $k$ tends to $+\infty$.\\
So there exists a sequence of integers $n_k$ tending to $+\infty$ and a sequence $\e_k$ tending to 0 such that :
$$\left\|(-i\partial_t+(-\Delta+V(x)+\kappa\varphi_k(x)))\left(e^{-i(\lambda_{n_k}+\kappa)t} u_{n_k}(x)\right)\right\|_{L^2(M)}
\leq \e_k+\kappa \left(1-\nu(\Gamma)+\e_k\right)^{\frac{1}{2}}$$
Put
$$v_k(t,x)=(-i\partial_t+(-\Delta+V(x)+\kappa\varphi_k(x)))\left(e^{-i(\lambda_{n_k}+\kappa)t} u_{n_k}(x)\right)$$
Duhamel's formula gives 
$$\left(e^{-it(-\Delta+V+\kappa\varphi_{n_k})} u_{n_k}\right)(t,x)=\left(e^{-i(\lambda_{n_k}+\kappa)t} u_{n_k}(x)\right)
-i\int_0^t e^{i(s-t)(-\Delta+V+\kappa\varphi_{n_k})}v_k(s,x)\,ds$$
So 
$$\left\|e^{-it(-\Delta+V+\kappa\varphi_{n_k})} u_{n_k}-\left(e^{-i(\lambda_{n_k}+\kappa)t} u_{n_k}\right)\right\|_{L^2(M)}
\leq t\left(\e_k+\kappa \left(1-\nu(\Gamma)+\e_k\right)^{\frac{1}{2}}\right)$$
Similarly
$$\left\|e^{-it(-\Delta+V)} u_{n_k}-e^{-i\lambda_{n_k}t} u_{n_k}\right\|_{L^2(M)}
\leq \e_k\,t$$
Obviously 
$$\left\|e^{-i(\lambda_{n_k}+\kappa)t} u_{n_k}-e^{-i\lambda_{n_k}t} u_{n_k}\right\|_{L^2(M)}=\left|e^{i\kappa t}-1\right|$$
Hence 
$$\left\|e^{-it(-\Delta+V+\kappa\varphi_{n_k})} u_{n_k}-e^{-it(-\Delta+V)} u_{n_k}\right\|_{L^2(M)}
\geq\left|e^{i\kappa t}-1\right|-\left(2\e_k+\kappa \left(1-\nu(\Gamma)+\e_k\right)^{\frac{1}{2}}\right)t$$
So 
$$\liminf_{k\rg +\infty }\left\|e^{-it(-\Delta+V+\kappa\varphi_{n_k})}-e^{-it(-\Delta+V)} \right\|_{\mathcal L(L^2(M))}
\geq 2\sin\frac{\kappa t}{2}-(1-\nu(\Gamma))^{\frac{1}{2}}\kappa\,t $$
For $\kappa$ small enough the right hand side is positive.\\

\begin{remark}  For $p=\infty$ there is always stability. Indeed, Duhamel's formula gives for $u_0$ in the domain of $-\Delta+V$
\ben e^{-it(-\Delta+V+W_n)}u_0-e^{-it(-\Delta+V)}u_0
=- i \int_0^t e^{-i(t-s)(-\Delta+V+W_n)} \, W_n \, e^{-is(-\Delta+V)}u_0\,ds \een
so 
\ben \|e^{-it(-\Delta+V+W_n)}-e^{-it(-\Delta+V)}\|_{\mathcal L (L^2(M))}\leq t\|W_n\|_{L^\infty(M)} \een
\end{remark}

\begin{remark} D. Jakobson (refeering to a communication of J. Bourgain) proves in \cite{Jak} that the second hypothetis of
 Theorem \ref{mainth} is not satisfied in the case of the Laplacian on tori. The stability is proved in dimension $2$ in \cite{BBZ},
the problem is open for $d\geq 3$.
\end{remark}

\section{Examples of applications} \label{exa}

 \subsection{The sphere $\mathbb S^d$}

We consider the sphere $\mathbb S^d$, $d\geq 2$ endowed with the usual metric induced by that one of $\R^{d+1}$.
We will show that Theorem \ref{mainth} gives the unstability near $-\Delta_{\mathbb S^d}$.\\
Let $e_n$ be the restriction to $\mathbb S^d$ of the harmonic
polynomial $(x_1+i\,x_2)^n$ (called equatorial spherical harmonic) and $u_n=e_n/\|e_n\|_{L^2(\mathbb S^d)}$
which satisfies
$$-\Delta_{\mathbb S^d}\,u_n=n(n+d-1)\,u_n,~~\|u_n\|_{L^2(\mathbb S^d)}=1$$
and we have to study the weak* limit of $|u_n|^2 d\mu$ where $d\mu$ is the Riemannian volume form on $\mathbb S^d$.

\begin{proposition} The sequence of measures $|u_n(x)|^2d\mu$ on $\mathbb S^d$ tends to the measure $\ds\frac{d\theta}{2\pi}$
on the circle $\{(\cos\theta,\sin\theta,0,\cdots,0)\,;\,\theta\in[0,\pi]\}$ for the weak* topology.
\end{proposition}

{\bf Proof.} We parameter $\mathbb S^d$ by
\ben x=(\cos\theta \cos\varphi,\sin\theta \cos\varphi,(\sin\varphi) \,t),~~
(\theta,\varphi,t)\in[0,2\pi]\times[0,\frac{\pi}{2}]\times \mathcal S^{d-2}
      \een
so  $d\mu=\cos\varphi\,(\sin\varphi)^{d-2}d\theta\,d\varphi\,dt$.\\
Let $f$ be a continuous fonction on $\mathbb S ^d$. 
\ben \int_{\mathcal S^d} f\,|e_n|^2\,d\mu=\int_0^{\frac{\pi}{2}} g(\varphi)
(\cos\varphi)^{2n+1}(\sin\varphi)^{d-2} d\varphi \een
where
\ben
g(\varphi)=\int_{[0,2\pi]\times \mathcal S^{d-2}} f(\cos\theta \cos\varphi,\sin\theta\cos\varphi,(\sin\varphi)\,t)\,d\theta\,dt
\een
Laplace's methods gives 
\ben
\int_{\mathbb S^d} f\,|e_n|^2\,d\mu=\left(\int_0^{+\infty}e^{-(n+\frac{1}{2})\varphi^2} \varphi^{d-2} d\varphi\right)
(g(0)+o(1))
\een
Notice that $\|e_n\|^2_{L^2(\mathbb S^d)}$ is given by $f=1$ so 
\ben \int_{\mathbb S^d} f(x)\,|u_n(x)|^2\,d\mu=\frac{1}{2\pi}\int_0^{2\pi} f(\cos\theta,\sin\theta,0,\cdots,0)\,d\theta+o(1)\een
which achieves the proof.\\

Thus Theorem \ref{mainth} applies on the sphere near $-\Delta_{\mathbb S^d}$.

\begin{remark} The result applies whenever the manifold $(M,g)$ only coincide with the sphere near a closed geodesic.
\end{remark}

\subsection{Periodic stable geodesic}

We refer here to the work of J.V.Ralston \cite{Ral3} who developed the ideas of \cite{Ba} and \cite{CdV} to obtain quasimodes using
WKB constructions. The author considers a Riemannian manifold $(M,g)$ of dimension $d\geq 2$ and assumes the existence of a
periodic closed non degenerate geodesic $\gamma$. The eigenvalues $\lambda_j$, $1\leq j\leq2(d-1)$ of the
Poincar\'e application assciated to $\gamma$ are supposed to have modulus $1$ (so $\gamma$ is stable) and to satisfy
the diophantian condition
\ben \forall n\in\N^{2d-2},~~\prod_{j=1}^{2d-2} \lambda_j^{n_j} \not=1 \een
$V$ is a smooth potential on $M$. Under these hypotheses, J.V. Ralston proves

\begin{theorem} For any nonnegative integer $N$ and real number $\varepsilon$, there exist sequences of
quasimodes $E_j$ tending to
$+\infty$ and associated normalized quasi-eigenfunctions $u_n$, and a constant $C_{N,\varepsilon}$ suth that
\begin{itemize}
 \item[1.]  $\|(-\Delta_g +V-E_n)\,u_n\|_{L^2(M)}=\mathcal O(E_n^{-N})$
\item[2.] $ \|u_n\|_{L^2(\gamma_{N,\varepsilon})}<\varepsilon $, where
$\gamma_{N,\varepsilon}=\{x\in M; d(x,\gamma)>C_{N,\varepsilon}E_n^{-\nicefrac14}\}$.
\end{itemize}
\end{theorem}

As a consequence, Theorem \ref{mainth} applies on $M$ for any smooth nonnegative potential $V$.

\subsection{Hyperbolic surfaces}

We refer here to the work of Y. Colin de Verdi\`ere and B. Parisse \cite{CdVP1}. They state their theorem in a particular case,
but it is easy to check that they proved the following 

\begin{theorem} \label{cdvpth}
Consider $a<b$ two real numbers, $f$ a smooth positive function on $[a,b]$, $H$ the cylinder $[a,b]\times\mathbb Z/2\pi\mathbb Z$
endowed with the metric $dt^2+f^2(t)\,d^2\theta$, and $\Delta$ the Laplace-Beltrami operator on $H$ associated
with this metric 
with Dirichlet boundary condition. Assume that $f$ reaches its minimum on $[a,b]$ in a unique point $t_0\in]a,b[$, and that 
this minimum is not degenerate.\\
Then $\Delta$ admits a sequence of normalized eigenfunctions $u_n$ suth that the sequence of measures $|u_n|^2d\mu$ tends to
$\ds\delta_{t_0}\otimes \frac{d\theta}{2\pi}$ for the weak* topology.
\end{theorem}

In particular, we can consider a rotation invariant surface. Let $a<b$ two real numbers, $f$ and $g$ two smooth functions
on $[a,b]$ such that $f>0$ and the curve $\mathcal C\,:\,[a,b]\ni t\mapsto (g(t),f(t),0)$ in $\mathbb R^3$ have no multiple
point and is
parametrized by its curvi-linear abscissa, i.e. $f'^2+g'^2=1$. The surface
$$H=\{(g(t),f(t)\cos\theta,f(t)\sin\theta);\,(t,\theta)\in[a,b]\times\mathbb Z/2\pi\mathbb Z\}$$
generated by the rotation of $\mathcal C$ on the X-axis is endowed with the metric induced by that one of $\mathbb \R^3$
$dt^2+f^2(t)\,d\theta^2$. Theorem \ref{mainth} applies to the Laplacian on $H$ with Dirichlet boundary condition.

\section{Unstability for small $p$} \label{ppetit}

In the case when $p<\dim M/2$, the unstability can be proved by local constructions and does not request the geometrical hypothetis.

\begin{theorem} \label{lowdth}
Let $(M,g)$ be a Riemannian manifold of dimension $d\geq 3$, $\Delta$ the associated Laplace-Betrami operator, $V$
a smooth nonnegative potential.\\
Then there exists a sequence of smooth bounded potentials $W_n$ such that
$$\forall p\in [1,\frac{d}{2}[,~~\lim_{n\rg +\infty} \|W_n\|_{L^{p}(M)}=0$$
$$\forall T>0,~~\lim_{n\rg +\infty} \|e^{-it(-\Delta+V+W_n)}- e^{-it(-\Delta+V)}\|_{L^\infty ([0,T];\,\mathcal L(L^2(M)))} =2$$
\end{theorem}

\begin{remark} Here the smooth potential $V$ is estimated in $L^p$-norm. Notice that the threshold $d/2$ is natural, 
 as it is the same one above which the operator $-\Delta+V$ in well defined (in the quadratic form sense).
\end{remark}

{\bf Proof of Theorem \ref{lowdth}:} In local coordinates we have
\be \Delta =\frac{1}{\sqrt{\det G}}\nabla \sqrt{\det G}\, G^{-1}\nabla  \label{defdelta} \ee
where $G$ is the matrix of the metric $g$.\\
Let $u^0\in C_0^\infty(\mathbb R^d)$ normalized by $\ds\|u^0\|_{\L^2(\mathbb R^d)}=1$ and
$W\in C_0^\infty(\mathbb R^d)$ real valued such that $W=1$ on the support of $u^0$.\\
We define for $n\in\mathbb N^*$ 
\ba u^0_n(\cdot)&=&n^{d/2} u_0(n\,\cdot) \label{defu} \\
W_n(\cdot)&=&n^2\ln (n+1)\, W_0(n\,\cdot) \label{defw} \ea
 so $\ds \|u^0_n\|_{L^2(\mathbb R^d)}=1$ and $\ds \|W_n\|_{L^2(\mathbb R^d)}=n^{2-\frac dp}\ln (n+1)$ tends to $0$ for any $p<d/2$.\\
By Duhamel's formula we have for any posive $t$
\be e^{-it(-\Delta+ V)} u^0_n-e^{-it V} u^0_n =
i\int_0^t e^{-i(t-s)(-\Delta+ V)} \Delta\left(e^{-is V} u^0_n\right)\,ds \label{du1} \ee
\be e^{-it(-\Delta+ V+W_n)} u^0_n-e^{-it(V+W_n)} u^0_n 
=i\int_0^t e^{-i(t-s)(-\Delta+V+W_n)} \Delta \left(e^{-is( V+W_n)} u^0_n\right)\,ds \label{du2} \ee
Definitions (\ref{defu}) and (\ref{defw}) of $u_n^0$ and $W_n$, together with (\ref{defdelta}) yield for a suitable $C>0$
\ba  \|\Delta\left(e^{-is V} u^0_n\right)\|_{L^2(\mathbb R^d)} &\leq& C\,n^2 \label{cc} \\
\|\Delta\left(e^{-is( V+W_n)} u^0_n\right)\|_{L^2(\mathbb R^d)} & \leq & C(s^2n^6\ln^2(n+1)+n^2) \label{dd} \ea
Then (\ref{du1}) and (\ref{du2}) give
\ba \|e^{-it(-\Delta+\ V)} u^0_n-e^{-it V} u^0_n\|_{L^2(\mathbb R^d)} &\leq& C\,n^2t \label{est1}\\
\|e^{-it(-\Delta+ V+W_n)} u^0_n-e^{-it( V+W_n)} u^0_n\|_{L^2(\mathbb R^d)} & \leq & C(t^3n^6\ln^2(n+1)+n^2t) \label{est2}
\ea
We consider the time 
\be t_n=\frac{\pi}{n^2\ln(n+1)} \ee
From (\ref{est1}) and (\ref{est2}) we get
\ba \|e^{-it_n(-\Delta+ V)} u^0_n-e^{-it_n V} u^0_n\|_{L^2(\mathbb R^d)} &\leq& \frac{C}{\ln(n+1)} \label{est1b}\\
\|e^{-it_n(-\Delta+ V+W_n)} u^0_n-e^{-it_n( V+W_n)} u^0_n\|_{L^2(\mathbb R^d)} & \leq & \frac{C}{\ln(n+1)} \label{est2b}
\ea
Since $W_n=1$ on the support of $u_n^0$ and $\ds\|u_n^0\|_{L^2(\mathbb R^d)}=1$ we have
\be \|e^{-it_n( V+W_n)} u^0_n - e^{-it_n V} u^0_n\|_{L^2(\mathbb R^d)}=2 \label {est3} \ee
So (\ref{est1}), (\ref{est2}) and (\ref{est3}) give
\be \|e^{-it_n(-\Delta+ V+W_n)} u^0_n - e^{-it_n(-\Delta+ V)} u^0_n\|_{L^2(\mathbb R^d)} \geq 2-\frac{C}{\ln(n+1)} \label{est4} \ee
which achieves the proof.

\nocite{CdVP2}

\bibliographystyle{plain}
\bibliography{unstability}

\end{document}